\documentclass[12pt]{article}

\newcommand{\blind}{1}

\usepackage{enumerate}
\usepackage{rotating}
\usepackage{orcidlink}
\usepackage{natbib}
\usepackage{pdflscape,cancel,latexsym,graphics,subcaption,color}
\usepackage{amsfonts,amsmath,amssymb,amsthm}
\usepackage{mathtools,xparse}
\usepackage{upgreek}
\usepackage[normalem]{ulem} 
\usepackage{pifont}
\RequirePackage{lineno}
\usepackage{hyperref}
\usepackage[margin=1in]{geometry}
\usepackage{longtable}
\usepackage{booktabs}
\usepackage{soul}

\usepackage{bm}

\newtheorem{proposition}{Proposition}
\newtheorem{lemma}{Lemma}
\newtheorem{corollary}{Corollary}

\newtheorem{definition}{Definition}
\newtheorem{example}{Example}

\definecolor{dred}{rgb}{0.5,0,0}
\definecolor{dgreen}{rgb}{0,0.7,0}
\definecolor{dblue}{rgb}{0,0,0.85}
\definecolor{dmagenta}{cmyk}{0,1,0,0.6}
\definecolor{dcyan}{cmyk}{1,0,0,0.5}
\definecolor{grey}{gray}{0.9}
\definecolor{orange}{rgb}{1,0.65,0}
\definecolor{mellow}{rgb}{.847,.72,.525}
\definecolor{golden}{rgb}{.80392,.60784,.11373}
\definecolor{dgolden}{rgb}{.5451,.39608,.03137}
\definecolor{brown}{rgb}{.15,.15,.15}
\definecolor{darkolivegreen}{rgb}{.33333,.41961,.18431}

\def\bone{{\bf 1}}

\def\cov{\mathop{\rm cov}\nolimits}
\def\tr{\mathop{\rm tr}\nolimits}

\newcommand{\dM}{\mathbb{M}}

\newcommand{\dO}{\mathbb{O}}
\newcommand{\dP}{\mathbb{P}}

\def\mR{\mathbb{R}}

\def\balpha{{\bm\alpha}}

\def\bvarepsilon{{\bm\varepsilon}}

\def\btheta{{\bm\theta}}

\def\blambda{{\bm\lambda}}

\def\bPsi{{\bm\Psi}}

\def\bff{{\bf f}}

\def\bp{{\bf p}}

\def\bw{{\bf w}}

\def\by{{\bf y}}

\def\bH{{\bf H}}
\def\bI{{\bf I}}

\def\bU{{\bf U}}
\def\bV{{\bf V}}
\def\bW{{\bf W}}

\def\cA{{\mathcal A}}
\def\cB{{\mathcal B}}
\def\cC{{\mathcal C}}

\def\cF{{\mathcal F}}
\def\cG{{\mathcal G}}
\def\cH{{\mathcal H}}

\def\cM{{\mathcal M}}
\def\cN{{\mathcal N}}

\def\cT{{\mathcal T}}

\def\cV{{\mathcal V}}

\def\cX{{\mathcal X}}

\DeclarePairedDelimiter{\abs}{\lvert}{\rvert}
\DeclarePairedDelimiter{\norm}{\lVert}{\rVert}

\DeclareMathOperator{\MVN}{MVN}
\def\ve{\varepsilon}

\newcommand{\fbm}{{FBM}}
\newcommand{\LL}{L} 
\newcommand{\N}{\operatorname{N}}
\newcommand{\E}{\operatorname{E}}
\newcommand*\diff{\mathop{}\!\mathrm{d}}

\usepackage{etoolbox} 
\makeatletter 
\newcommand*\linenomathpatch{\@ifstar{\linenomathpatch@AMS}{\linenomathpatch@}}
\newcommand*\linenomathpatch@[1]{
	\expandafter\pretocmd\csname #1\endcsname {\linenomathWithnumbers}{}{}
	\expandafter\pretocmd\csname #1*\endcsname{\linenomathWithnumbers}{}{}
	\expandafter\apptocmd\csname end#1\endcsname {\endlinenomath}{}{}
	\expandafter\apptocmd\csname end#1*\endcsname{\endlinenomath}{}{}
}
\newcommand*\linenomathpatch@AMS[1]{
	\expandafter\pretocmd\csname #1\endcsname {\linenomathWithnumbersAMS}{}{}
	\expandafter\pretocmd\csname #1*\endcsname{\linenomathWithnumbersAMS}{}{}
	\expandafter\apptocmd\csname end#1\endcsname {\endlinenomath}{}{}
	\expandafter\apptocmd\csname end#1*\endcsname{\endlinenomath}{}{}
}
\let\linenomathWithnumbersAMS\linenomathWithnumbers
\patchcmd\linenomathWithnumbersAMS{\advance\postdisplaypenalty\linenopenalty}{}{}{}
\makeatother 

\linenomathpatch{equation}
\linenomathpatch*{gather}
\linenomathpatch*{multline}
\linenomathpatch*{align}
\linenomathpatch*{alignat}
\linenomathpatch*{flalign}

\begin{document}

\def\spacingset#1{\renewcommand{\baselinestretch}%
{#1}\small\normalsize} \spacingset{1}


\if1\blind
{
  \title{\bf Additive interaction modelling using I-priors}
  \author{Wicher Bergsma \\
    London School of Economics and Political Science\\
    and \\
    Haziq Jamil \\
    Universiti Brunei Darussalam}
  \maketitle
} \fi

\if0\blind
{
  \bigskip
  \bigskip
  \bigskip
  \begin{center}
    {\LARGE\bf Additive interaction modelling using I-priors}
\end{center}
  \medskip
} \fi

\bigskip
\begin{abstract}
Additive regression models with interactions are widely studied in the literature, using methods such as splines or Gaussian process regression. However, these methods can pose challenges for estimation and model selection, due to the presence of many smoothing parameters and the lack of suitable criteria. We propose to address these challenges by extending the I-prior methodology \citep{bergsma20} to multiple covariates, which may be multidimensional. 
The I-prior methodology has some advantages over other methods, such as Gaussian process regression and Tikhonov regularization, both theoretically and practically. In particular, the I-prior is a proper prior, is based on minimal assumptions, yields an admissible posterior mean, and estimation of the scale (or smoothing) parameters can be done using an EM algorithm with simple E and M steps. Moreover, we introduce a parsimonious specification of models with interactions, which has two benefits: (i) it reduces the number of scale parameters and thus facilitates the estimation of models with interactions, and (ii) it enables straightforward model selection (among models with different interactions) based on the marginal likelihood.
\end{abstract}

\noindent%
\vfill

\newpage

\section{Introduction}\label{sec-intro}


In this paper we consider a novel methodology for estimating and selecting regression models with interactions. 
For an observed covariate $x$, a scalar response $y$, and an error term $\ve$, interaction models have the form $y=f(x)+\ve$, where $f$ satisfies an additive structure with some possible interactions between covariate components.
As an example, if $x=(x_1,x_2,x_3)$ has three components, the {\em saturated} model, with all interaction effects included, is
\begin{align} \label{sat3}
	\begin{split}
	  f(x) &= \alpha + \underbrace{f_1(x_{1}) + f_2(x_{2}) + f_3(x_{3})}_{\text{main effects}} + \underbrace{f_{12}(x_{1},x_{2}) + f_{23}(x_{2},x_{3}) + f_{13}(x_{1},x_{3})}_{\text{1st order interactions}} \\
    &\quad + \underbrace{f_{123}(x_{1},x_{2},x_{3})}_{\text{2nd order interaction}}.
	\end{split}
\end{align}
Here, the $x_{k}$ are possibly multidimensional covariates, $\alpha$ is a real number, and the subscripted $f$ are real-valued functions. 
Of interest is to know whether some of the interactions are zero, e.g., we may wish to know whether the regression function is instead of the following form:
\begin{align} \label{nonsat}
	  f(x) = \alpha + f_1(x_{1}) + f_2(x_{2}) + f_3(x_{3}) + f_{12}(x_{1},x_{2}) + f_{13}(x_{1},x_{3}).
\end{align}
In this model, the effects of $x_2$ and $x_3$ on $y$ both depend on $x_1$, but do not depend on $x_3$ and $x_2$, respectively. 
In practice, like~(\ref{nonsat}), we require models to be {\em hierarchical}, that is, if a higher order interaction is present then so are the corresponding lower order interactions.

The initial approach to additive modelling was to approximate the regression function using spline functions, obtained by minimizing the sum of squared errors subject to an appropriate penalty (e.g., \cite{stone85,wahba86,wahba90,friedman91,gu13}). 
The spline methodology is typically based on quadratic form penalties and is then also referred to as Tikhonov regularization. A well-known (empirical or fully) Bayesian interpretation of spline models exists \citep{kw70bayes}. In recent years, Gaussian process regression is also gaining popularity for estimating additive models with interactions (e.g., \cite{rw06}). This approach involves assuming a Gaussian process prior for the regression function, and estimating the regression functions by their posterior mean. In view of the Bayesian interpretation of spline smoothers, the latter can be viewed as posterior means under a Gaussian process prior.


In this paper we introduce an alternative approach to additive interaction modelling which has some advantages compared to the aforementioned ones, in particular, we extend the I-prior methodology introduced by \citet{bergsma20} to be able to cover multiple covariates. \citet{jamil18} provides a number of extensions to the methodology introduced here, including probit and logit models using a fully Bayes approach and Nystr\"om approximations for speeding up the I-prior methodology. Furthermore, he contributed a user friendly R package \texttt{iprior} \citep{jamil19package}, further described in \citet{jb19}.

The I-prior methodology is a {\em function space approach}, that is, the regression function is assumed $f$ to lie in a particular space of functions. 
The I-prior maximizes entropy, subject to a constraint based on the Fisher information on $f$. Its support is contained in the function space, i.e., the I-prior is proper. If the  errors are normal, the I-prior is Gaussian with covariance kernel proportional to the Fisher information on the regression function. I-priors have the intuitively appealing property that the more information is available on a linear functional of the regression function, the larger the prior variance, and the smaller the influence of the prior mean on the posterior distribution.

The space we assume the regression function to be in is a {\em reproducing kernel Krein space} (RKKS). These spaces generalize reproducing kernel Hilbert spaces (RKHSs), dropping the positive definiteness restriction on the reproducing kernel. The Fisher information of the regression function is positive definite regardless of the definiteness of the kernel, hence, there is no need to restrict attention to RKHSs. In fact, this would be arbitrarily and unnecessarily restrictive for I-prior modelling with multiple covariates. In this case, the kernel is a weighted sum of (usually positive definite) kernels corresponding to each covariate; there is no reason to restrict the weights to be positive. 

The main competing methods to the I-prior methodology are Gaussian process (GP) regression and Tikhonov regularization and the I-prior has some theoretical and practical advantages compared to these. Firstly, in contrast to I-prior estimators, Tikhonov regularizers are inadmissible with respect to squared error loss \citep{cp19}. Secondly, in GP regression, both the space in which the regression function lives and a prior over that space needs to be specified; only the former is needed for I-priors, which are then generated automatically. Practical advantages of our methodology are ease of estimation and model comparison.
In particular, we develop an EM algorithm with a simple E and M step for estimating certain hyperparameters, facilitating estimation for complex models. We also propose a novel parsimonious model formulation, requiring a single scale parameter for each (possibly multidimensional) covariate and no further parameters for interaction effects. This simplifies estimation because fewer hyperparameters need to be estimated, and also simplifies model comparison of models with the same covariates but different interaction effects; in this case, the model with the highest estimated likelihood can be selected.

I-priors are equally suitable for parametric and nonparametric regressions. We show via simulations that, with univariate covariates and linear effects, that is, for saturated model $y=f(x)+\ve$ where
\begin{align}\label{sat3-lin}
	f(x) = \alpha + x_{1}\beta_{1} + x_{2}\beta_{2} + x_{3}\beta_{3} + x_{1}x_{2}\beta_{12} + x_{2}x_{3}\beta_{23} + x_{1}x_{3}\beta_{13} + x_{1}x_{2}x_{3}\beta_{123} ,
\end{align}
and in the settings we investigated, the I-prior methodology has typically higher probability of selecting the correct hierarchical model than competing methods such as lasso, spike-and-slab, or $g$-priors. However, the I-prior methodology can also straightforwardly be used in more general settings, which the competing methods cannot.





The present paper extends \citet{bergsma20}, which covers the case of a single, possibly multidimensional covariate. 
In this case there is a single scale parameter for which the RKHS framework suffices, whereas in the present paper we require possibly indefinite RKKSs.
The latter paper outlines the relation with competing methods, including $g$-priors (which can be viewed as I-priors for linear regression functions when the covariate space is equipped with Mahalanobis distance), Jeffreys and reference priors, and Fisher kernels. 
A detailed comparison with Tikhonov regularization is given, with particular detail on the relation with cubic spline smoothing. 
It is explained how I-priors work when the regression functions are linear, or when they are assumed to lie in the fractional Brownian motion RKHS.

\citet{ong04} used RKKSs in a regularization framework, where the usual RKHS squared penalty norm $\norm{f}_\cF^2$ is replaced by the RKKS indefinite inner product $\langle f,f\rangle_\cF$. As the latter may be negative, it does not make sense to minimize the ``penalized'' loss function, and instead they sought a saddle point. In contrast, in the I-prior approach a Gaussian prior for the regression function whose covariance kernel is proportional to its Fisher information is used; the Fisher information is positive definite regardless of the definiteness of the reproducing kernel. 


This paper is organised as follows.
Section~\ref{sec-rkks} gives a brief overview of RKKSs and we provide lemmas on sums and tensor products of RKKSs, as needed for Section~\ref{sec-add}. There, hierarchical interaction spaces are defined formally and in some generality. Based on this, our definition of hierarchical interaction models is then given.
In Section~\ref{sec-noninf}, the I-prior is defined for the regression function in a hierarchical interaction model when the regression function lies in an RKKS.
An efficient EM algorithm for estimating scale parameters is described. 
In Section~\ref{sec-sim}, simulations are done to evaluate model selection performance in the simple setting of linear regression functions with interactions. These show good performance compared to selection using the lasso, spike \& slab priors or $g$-priors. 
In Section~\ref{sec-cow}, we apply the I-prior methodology to a real data example, illustrating model selection when the regression function is nonlinear. 
The paper concludes with a discussion in Section~\ref{sec-disc}. 

\section{Reproducing kernel Krein spaces}\label{sec-rkks}

{\em Krein spaces} generalize Hilbert spaces by relaxing the positive definiteness restriction on the inner product, while
{\em Reproducing kernel Krein spaces} (RKKSs) generalize {\em reproducing kernel Hilbert spaces} by relaxing the positive definiteness restriction on the kernel. An RKKS possesses a unique symmetric reproducing kernel (r.k.) which is the difference of two positive definite kernels. 
Below, we give a brief overview of the theory as needed for this paper, giving proofs if we could not find a reference for a statement. 
Further details can be found in \citet{ong04,bognar74,azizov89,gheondea13}.

\begin{definition}[Inner product]
Let $\cF$ be a vector space over the reals.
A function $\langle\cdot,\cdot,\rangle_\cF:\cF\times\cF\rightarrow\mR$ is called an {\em inner product} on $\cF$ if, for all $f,f',f''\in\cF$,
\begin{itemize}

\item (symmetry) $\langle f,f'\rangle_\cF=\langle f',f\rangle_\cF$

\item (linearity) $\langle \alpha f+f',f''\rangle_\cF=\alpha\langle f,f''\rangle_\cF+\langle f',f''\rangle_\cF$

\item (nondegeneracy) $(\forall g\in\cF:\langle f,g\rangle_\cF=0)\Rightarrow f=0$

\end{itemize}
If $\langle f,f\rangle_\cF> 0$ for all $f\in\cF$, $f\ne 0$, the inner product is called {\em positive definite}.
\end{definition}

Recall that a Hilbert space is a complete positive definite inner product space.
The more general notion of Krein space is defined as follows.
\begin{definition}[Krein space]\label{def-krein}
A vector space $\cF$ equipped with the inner product $\langle\cdot,\cdot,\rangle_\cF$ is called a {\em Krein space} if there are two Hilbert spaces $\cH_+$ and $\cH_-$ with $\cH_+\cap\cH_-=\{0\}$ spanning $\cF$ such that
\begin{itemize}
\item all $f\in\cF$ can be decomposed as $f=f_+ +f_-$ where $f_+\in\cH_+$ and $f_-\in\cH_-$.
\item for all $f,f'\in\cF$, $\langle f,f'\rangle_\cF=\langle f_+,f_+'\rangle_{\cH_+} - \langle f_-,f_-'\rangle_{\cH_-}$.
\end{itemize}
\end{definition}
\noindent Thus, any Hilbert space is a Krein space, which can be seen by taking $\cH_-=\{0\}$. Note that, for any Krein space, there exists a Hilbert space containing the same functions:
\begin{definition}[Associated Hilbert space]
Let $\cF$ be a Krein space with decomposition into Hilbert spaces $\cH_+$ and $\cH_-$. and inner product $\langle\cdot,\cdot\rangle_{\cH_+}-\langle\cdot,\cdot\rangle_{\cH_-}$. Then we denote by $\overline\cF$ the associated Hilbert space $\cH_+\oplus\cH_-$ with the inner product $\langle\cdot,\cdot\rangle_{\cH_+}+\langle\cdot,\cdot\rangle_{\cH_-}$.
\end{definition}
There may be multiple Hilbert spaces associated to a Krein space. The {\em strong topology} of a Krein space $\cF$ is defined as the Hilbert topology of $\overline\cF$, and does not depend on the chosen decomposition of $\cF$.
The following definition is based on \citet{gheondea13}.
\begin{definition}[RKKS]
A Krein space $\cF$ consisting of functions $f:\cX\to\mR$ is called a {\em reproducing kernel Krein space} (RKKS) if there exists a symmetric function $h:\cX\times\cX\to\mR$, such that, for all $x\in\cX$,
\begin{itemize}
    \item $h(x,\cdot)\in\cF$, and
    \item for all $f\in\cF$, $f(x)=\langle f,h(x,\cdot)\rangle_{\cF}$
\end{itemize}
The function $h$ is called the {\em reproducing kernel} (r.k.) of $\cF$.
\end{definition}
Every RKKS has a unique r.k.\ which is a difference of two positive definite kernels \citep[Proposition 6]{ong04}, but certain kernels may have multiple RKKS associated with them \citep{alpay91,gheondea13}. 

We will now show how sums of disjoint RKKSs and tensor products of arbitrary RKKSs can be constructed. We use the following notation: for a Hilbert space $\cH$ with inner product $\langle\cdot,\cdot\rangle_{\cH}$, $\cH^-$ denotes the Krein space with inner product $-\langle\cdot,\cdot\rangle_{\cH}$. 



\begin{lemma}\label{lem-decomp}
    Let $\cH_+$ and $\cH_-$ be RKHSs on $\cX$ with r.k.s $h_+$ and $h_-$ such that $\cH_+\cap\cH_-=\{0\}$. 
    Then the space $\cF=\cH_+\oplus\cH_-^-=\{f|f=f_++f_-, f_+\in\cH_+, f_-\in\cH_-\}$ equipped with the inner product
    \[ \langle f,f'\rangle_{\cF} = \langle f_+,f_+' \rangle_{\cH_+}-\langle f_-,f_-'\rangle_{\cH_-} \]
    is an RKKS with r.k.\ $h_+-h_-$.
\end{lemma}
\begin{proof}
    By definition, $\cF$ is a Krein space. Clearly, $(h_+-h_-)(x,\cdot)=h_+(x,\cdot)-h_-(x,\cdot)\in\cF$. 
    Furthermore,
    \begin{align*}
        \langle f, h_+(x,\cdot) - h_-(x,\cdot) \rangle_\cF 
        =
        \langle f_+,h_+(x,\cdot) \rangle_{\cH_+} + \langle f_-,h_-(x,\cdot) \rangle_{\cH_-} = f_+(x)+f_-(x) = f(x).
    \end{align*}
    Hence, $h_+-h_-$ is an r.k.\ for $\cF$, completing the proof.
\end{proof}

\begin{proposition}[Sum of disjoint RKKSs] \label{prop-sum}
    Let $\cF_1$ and $\cF_2$ be RKKSs on $\cX$ with r.k.s $h_1$ and $h_2$ such that $\cF_1\cap\cF_2=\{0\}$. Then $\cF=\cF_1\oplus\cF_2=\{f|f=f_1+f_2,f_1\in\cF_1,f_2\in\cF_2\}$ equipped with the inner product
    \[ \langle f,f'\rangle_{\cF} = \langle f_1,f_1' \rangle_{\cF_1}+\langle f_2,f_2'\rangle_{\cF_2} \]
    is an RKKS with r.k. $h=h_1+h_2$.
\end{proposition}
\begin{proof}
    We assume notation for decompositions of Krein spaces analogous to the above. Then
    \begin{align*}
        \cF = \cF_1\oplus\cF_2 = (\cH_{1+}\oplus\cH_{1-}^-)\oplus(\cH_{2+}\oplus\cH_{2-}^-) = (\cH_{1+}\oplus\cH_{2+})\oplus (\cH_{1-}\oplus\cH_{2-})^- .
    \end{align*}
    By Theorem 5 of \citet{bt04}, $\cH_{1+}\oplus\cH_{2+}$ and $\cH_{1-}\oplus\cH_{2-}$ are RKHSs with r.k.s $h_{1+}+h_{2+}$ and $h_{1-}+h_{2-}$; their intersection is $\{0\}$. Hence, by Lemma~\ref{lem-decomp}, $\cF$ is an RKKS with kernel $h_{1+}+h_{2+}+h_{1-}+h_{2-}=h_1+h_2$.    
\end{proof}
As a remark, sums of non-disjoint RKKSs are not needed for this paper. 

For vector spaces $\cH_1$ and $\cH_2$, the tensor product $\cH_1\tilde\otimes\cH_2$ is defined as the vector space consisting of weighted sums of simple tensors $f_1\otimes f_2$ ($f_1\in\cH_1$, $f_2\in\cH_2$).
If $\cH_1$ and $\cH_2$ are Hilbert spaces, the norm of the simple tensor $f_1\otimes f_2$ ($f_1\in\cH_1$, $f_2\in\cH_2$) is defined as ${\norm{f_1\otimes f_2}_{\cH_1\otimes\cH_2}}=\norm{f_1}_{\cH_1}\norm{f_2}_{\cH_2}$.
The completion of $\cH_1\tilde\otimes\cH_2$ with respect to this norm is then denoted $\cH_1\otimes\cH_2$.
If $\cF_1$ and $\cF_2$ are RKHSs on $\cX_1$ and $\cX_2$ with r.k.s $h_1$ and $h_2$, then $\cF=\cF_1\otimes\cF_2$ can be shown to have r.k.\ $h=h_1\otimes h_2$ (\cite{bt04}, Theorem 13; \cite{neveu68}). Note that point evaluation of $h$ is given by
\begin{align*}
    h\big((x_1,x_2),(x_1',x_2')\big) = h_1(x_1,x_1')h_2(x_2,x_2') \quad \forall x_1,x_1'\in\cX_1,x_2,x_2'\in\cX_2.
\end{align*}

\begin{proposition}[Tensor products of RKKSs] \label{prop-tensor}
    Let $\cF_1$ and $\cF_2$ be RKKSs on on $\cX_1$ and $\cX_2$ with r.k.s $h_1$ and $h_2$. Then the tensor product $\cF_1\tilde\otimes\cF_2$ admits a functional completion $\cF_1\otimes\cF_2$ which is an RKKS with r.k.\ $h_1\otimes h_2$. 
\end{proposition}
\begin{proof}
With the notation for decompositions of Krein spaces as above, we have 
    \begin{align*}
        \cF_1\tilde\otimes\cF_2 
        &= (\cH_{1+}\oplus\cH_{1-}^-)\tilde\otimes(\cH_{2+}\oplus\cH_{2-}^-) \\
        &= (\cH_{1+}\tilde\otimes\cH_{2+}) \oplus (\cH_{1+}\tilde\otimes\cH_{2-}^-) \oplus (\cH_{1-}\tilde\otimes\cH_{2+}^-) \oplus (\cH_{1-}\tilde\otimes\cH_{2-}) \\
        &= \big((\cH_{1+}\tilde\otimes\cH_{2+}) \oplus (\cH_{1-}\tilde\otimes\cH_{2-})\big) \oplus \big( (\cH_{1+}\tilde\otimes\cH_{2-}) \oplus (\cH_{1-}\tilde\otimes\cH_{2+})\big)^-.
    \end{align*}
We may now define the functional completion of $\cF_1\tilde\otimes\cF_2$ as
\begin{align*}
    \cF_1\otimes\cF_2 
    &= \big(\underbrace{(\cH_{1+}\otimes\cH_{2+}) \oplus (\cH_{1-}\otimes\cH_{2-})}_{=:\cH_{12,+}}\big) \oplus \big( \underbrace{(\cH_{1+}\otimes\cH_{2-}) \oplus (\cH_{1-}\otimes\cH_{2+})}_{=:\cH_{12,-}} \big)^- .
\end{align*}
By Theorems~5 and~13 of \citet{bt04}, $\cH_{12,+}$ and $\cH_{12,-}$ are RKHSs with r.k.s $h_{12,+}:=h_{1+}\otimes h_{2+}+h_{1-}\otimes h_{2-}$ and $h_{12,-}:=h_{1+}\otimes h_{2-}+h_{1-}\otimes h_{2+}$, and $\cH_{12,+}\cap\cH_{12,-}=\{0\}$. 
Thus, by Lemma~\ref{lem-decomp}, $\cF_1\otimes\cF_2$ is an RKKS with kernel $h_{12,+} - h_{12,-} = (h_{1+}-h_{1-})\otimes(h_{2+}-h_{2-})$.
\end{proof}

\section{Hierarchical interaction models}\label{sec-add}

This section formalizes the concept of hierarchical interaction models.
Hierarchical interaction spaces are introduced in Section~\ref{sec-hier1} in some generality. To make them amenable to computation, Section~\ref{sec-kern} considers these spaces as RKKSs, and we derive the corresponding reproducing kernel. Our construction leads to a more parsimonious assignment of scale parameters than is typical in the literature. 
In Section~\ref{sec-hier2} we define hierarchical interaction models, which are the object of study of this paper. 

\subsection{Hierarchical interaction spaces}\label{sec-hier1}

Let $\cV$ be a finite set. We will refer to $\cV$ as a set of variables. The {\em power set} of $\cV$ is the set of all its subsets and is denoted $\dP(\cV)$. A subset $\dO\subseteq\dP(\cV)$ with the property that if $\cM\in\dO$ then also $\cM'\in\dO$ whenever $\cM'\subset\cM$ is called an {\em order ideal} of $\dP(\cV)$ (cf.\ \cite{stanley11}). 
The set of maximal elements of an order ideal of $\cV$ is called a {\em Sperner family} of $\cV$. 
A Sperner family is characterized by the property that it is a set of subsets such that none of its elements contains another. 
For $\abs{\cV}=0,1,2,3,4,5,\ldots$, the number of Sperner families of $\cV$ is given by the Dedekind numbers $2, 3, 6, 20, 168, 7581,\ldots$.  
For $\cV=\{1,2\}$, the 6 Sperner families and their corresponding order ideals are listed in Table~\ref{tbl-hier6}.

For $\dM\subseteq\dP(\cV)$, define $\tilde\dP(\dM)=\cup_{\cM\in\dM}\dP(\cM)$ where we understand $\tilde\dP(\emptyset)=\emptyset$.
Then if $\dM$ is a Sperner family, the corresponding order ideal is $\dO=\tilde\dP(\dM)$.
We shall define our models in terms of Sperner families rather than order ideals, because, having generally fewer elements, these are easier to represent.

For $v\in\cV$ let $\cX_v$ be a set, let $\cF_v$ be a vector space of real-valued functions on $\cX_v$, and let $\cC_v$ be the vector space of constant functions on $\cX_v$. We assume $\cF_v$ does not contain the constant functions, i.e., $\cF_v\cap\cC_v=\{0\}$. 
For $\cM\subseteq\cV$ and $v\in\cV$, let $t_\cM(v)=\cF_v$ if $v\in\cM$ and $t_\cM(v) =\cC_v$ if $v\in\cV\setminus\cM$.
Let $\cT_\cM = \otimes_{v\in\cV}t_\cM(v)$, where the tensor product is taken in the lexicographical order of the elements of $\cV$. For example, if $\cV=\{1,2\}$,
\[  \cT_{\emptyset} = \cC_1\otimes\cC_2, \quad 
    \cT_{\{1\}} = \cF_1\otimes\cC_2,\quad 
    \cT_{\{2\}} = \cC_1\otimes\cF_2, \quad \mbox{and }
    \cT_{\{1,2\}} = \cF_1\otimes\cF_2 . \]
\begin{definition}\label{def-his}
Let $\dM$ be a Sperner family of a finite set $\cV$ and for $v\in\cV$ let $\cF_v$ be a vector space of functions on a set $\cX_v$. Then the {\em hierarchical interaction space} $\cF_\dM$ is defined as
\[  \cF_\dM = \bigoplus_{\cM\in\tilde\dP(\dM)}\cT_\cM,
    \]
with the convention $\oplus_{\cM\in\emptyset}\cT_\cM=\{0\}$.
The {\em saturated interaction space} is defined as $\cF_{\{\cV\}}$.
\end{definition}

The following special case is of some interest:
\begin{corollary}
The saturated interaction space $\cF_{\{\cV\}}$ can be decomposed as 
\[  \cF_{\{\cV\}} = \bigotimes_{v\in\cV}\big(\cF_v\oplus\cC_v\big) . \]
\end{corollary}
\begin{proof}
Expanding the tensor product immediately gives
$\otimes_{v\in\cV}\big(\cF_v\oplus\cC_v\big) = \oplus_{\cM\subseteq\cV}\cT_\cM$.
The proof is completed by noting that $\{\cM \mid \cM\subseteq\cV\}=\dP(\cV)=\tilde\dP(\{\cV\})$.
\end{proof}

For non-empty $\cM\subseteq\cV$, denote $\cX_\cM=\times_{v\in\cM}\cX_v$, $\cF_\cM=\otimes_{v\in\cM}\cF_v$. For 
$x\in\cX_\cV$ and $\cM\subset\cV$, denote $x_\cM=(x_v)_{v\in\cM}$, that is, the tuple containing the coordinates of $x$ pertaining to $\cM$ (note $x_\emptyset=()$).
Then $f\in\cF_\dM$ can be represented in the form
\[ f(x) = \sum_{\cM\in\tilde\dP(\dM)} f_\cM(x_\cM) \quad\quad \forall x\in\cX_\cV,\]
where $f_\cM\in\cF_\cM$ and we use the convention that $\sum_{\cM\in\emptyset}f_\cM(x_\cM)=0$. 
For example, for the Sperner family  $\dM=\{\{1,2\},\{2,3\}\}$ of $\cV=\{1,2,3\}$, every $f\in\cF_\dM$ has the decomposition \vspace{-0.5em}
\[  f(x) = f_\emptyset() + f_{\{1\}}(x_{1}) + f_{\{2\}}(x_{2}) + f_{\{3\}}(x_{3}) + f_{\{1,2\}}(x_{1},x_{2}) + f_{\{2,3\}}(x_{2},x_{3}) \quad\forall x\in\cX_{\{1,2,3\}}.  \vspace{-0.5em} \]
In the sequel, we will write for convenience $f_{12}$ as shorthand for $f_{\{1,2\}}$, and so on. 
Table~\ref{tbl-hier6} gives the six hierarchical interaction spaces for $\cV=\{1,2\}$.


\begin{table}[t]
    \centering
    \caption{List of all Sperner families of $\cV=\{1,2\}$ and their corresponding power sets and hierarchical interaction spaces with the form of functions these contain.}    
    \label{tbl-hier6}
    \begin{tabular}{llll}
    \toprule
        Sperner family $\dM$           & Order ideal $\tilde\dP(\dM)$                    & $\cF_\dM$ / form of $f\in\cF_\dM$   \\ 
    \midrule
        $\emptyset$     & $\emptyset$                   & $\{0\}$ \\ & & $f(x_1,x_2)=0$  & \\
        $\{\emptyset\}$ & $\{\emptyset\}$               & $\cC_1\otimes\cC_2$ \\ & & $f(x_1,x_2)=f_\emptyset()$ &  \\
        $\{\{1\}\}$     & $\{\emptyset,\{1\}\}$         & $(\cC_1\otimes\cC_2) \oplus (\cF_1\otimes\cC_2$) \\ & & $f(x_1,x_2)=f_\emptyset()+f_1(x_1)$ & \\
        $\{\{2\}\}$     & $\{\emptyset,\{2\}\}$         & $(\cC_1\otimes\cC_2) \oplus (\cC_1\otimes\cF_2$) \\ & & $f(x_1,x_2)=f_\emptyset()+f_2(x_2)$ & \\
        $\{\{1\},\{2\}\}$ & $\{\emptyset,\{1\},\{2\}\}$ & $(\cC_1\otimes\cC_2) \oplus (\cF_1\otimes\cC_2)\oplus (\cC_1\otimes\cF_2)$ \\ & & $f(x_1,x_2)=f_\emptyset()+f_1(x_1)+f_2(x_2)$ & \\
        $\{\{1,2\}\}$   & $\{\emptyset,\{1\},\{2\},\{1,2\}\}$ & $(\cC_1\otimes\cC_2) \oplus (\cF_1\otimes\cC_2) \oplus (\cC_1\otimes\cF_2) \oplus (\cF_1\otimes\cF_2)$ \\ & & $f(x_1,x_2)=f_\emptyset()+f_1(x_1)+f_2(x_2)+f_{12}(x_1,x_2)$ & 
        \\ 
    \bottomrule
    \end{tabular}
\end{table}

\subsection{Kernel of hierarchical interaction spaces}\label{sec-kern}

Suppose that for $v\in\cV$, $\cF_v$ and $\cV_v$ are RKKSs with r.k.s $h_v$ and $c_v$.  
For $v\in\cV$ and $\cM\subseteq\cV$, let $u_\cM(v)=h_v$ if $v\in\cM$ and $u_\cM(v)=c_v$ if $v\in\cV\setminus\cM$, and let $k_\cM=\otimes_{v\in\cV}u_\cM(v)$. For example, if $\cV=\{1,2\}$,
 \[  k_{\emptyset} = c_1\otimes c_2, \quad 
    k_{\{1\}} = h_1\otimes c_2,\quad 
    k_{\{2\}} = c_1\otimes h_2, \quad \mbox{and }
    k_{\{1,2\}} = h_1\otimes h_2. \]
In the next proposition considers the kernel of the hierarchical interaction space $\cF_\dM$ (Definition~\ref{def-his}) when the components are RKKSs. We assume the direct sum and tensor products of RKKSs are as implied by  Propositions~\ref{prop-sum} and~\ref{prop-tensor}.

\begin{proposition}[Kernel of hierarchical interaction RKKS] \label{prop-kern}
If $\dM$ is a Sperner family of $\cV$ and $\cF_v$ and $\cC_v$ are RKKSs with r.k.s $h_v$ and $c_v$ such that $\cF_v\cap\cC_v=\{0\}$, then $\cF_\dM$ is the RKKS with kernel 
\[  h_\dM(x,x') = 
        \sum_{\cM\in\tilde\dP(\dM)}k_\cM(x_\cM,x_\cM').
\]
\end{proposition}
\begin{proof}
By Lemma~\ref{prop-tensor}, for $\cM\subseteq\cV$, $\cT_\cM$ is an RKKS with r.k.\ $k_\cM$. Any pair $\cT_\cM$ and $\cT_{\cM'}$ with $\cM\ne\cM'$ differ in at least one component, that is, there is a $v\in\cV$ such that the $v$th component of one is $\cF_v$ and $\cC_v$ for the other. Since the latter two are disjoint, the same is true for $\cT_{\cM}$ and $\cT_{\cM'}$. Hence, by Proposition~\ref{prop-sum}, their kernels can be summed in the required way, completing the proof. 
\end{proof}

In this paper we make the following assumptions for the kernels:
\begin{itemize}
\item[A1] For $v\in\cV$, $\cF_v$ is an RKKS on $\cX_v$ with r.k.\  $h_v=\tilde\lambda_v\tilde h_v$, where $\tilde h_v$ is a positive definite kernel and $\tilde\lambda_v$ is a real-valued parameter. Furthermore, $\cF_v$ does not contain any nonzero constant functions. 

\item[A2] For $v\in\cV$, $\cC_v$ is an RKKS of constant functions on $\cX_v$ with kernel $c_v(x,x')=\xi_v$ where $\xi_v\ne 0$ is a real-valued parameter.
\end{itemize}
We call the $\tilde\lambda_v$ and $\xi_v$ {\em scale parameters} as they determine the length scales of vectors in the RKHS with kernels $\tilde h_v$ and $\abs{c_v}$. In general these length scales are essentially arbitrary and not comparable across function spaces, hence these need to be estimated. Length scales in tensor product spaces inherit the length scales of the components through multiplication.

Under assumptions A1 and A2, the kernel of Proposition~\ref{prop-kern} evaluates to
\begin{align}\label{parsim kern unident}  
    h_\dM(x,x') = \sum_{\cM\in\tilde\dP(\dM)}\Big(\prod_{v\in{\cV\setminus\cM}}\xi_v\Big)\prod_{v\in\cM}\tilde\lambda_v\tilde h_v(x_v,x_v'),
\end{align}
where we use the convention $\prod_{v\in\emptyset}=1$. 
The scale parameters are not identified and we can identify them as follows. Set $\tau:=\prod_{v\in\cV}\xi_v$ and $\lambda_v:=\tilde\lambda_v/\xi_v$. Then the kernel reduces to 
\begin{align}\label{parsim kern}  
    h_\dM(x,x') = \tau\sum_{\cM\in\tilde\dP(\dM)}\prod_{v\in\cM}\lambda_v\tilde h_v(x_v,x_v').
\end{align}
It can easily be seen that identification is completed (i.e., different values of the scale parameters lead to different kernels) by assuming $\tau>0$. 
The saturated kernel then reduces to 
\begin{align}\label{parsim sat kern}  
    h_\dM(x,x') = \tau\prod_{v\in\cV}\big(1+\lambda_v\tilde h_v(x_v,x_v')\big) .
\end{align}
Thus, for any hierarchical interaction model which includes $\abs{\cV}$ main effects, the required number of scale parameters is $\abs{\cV}+1$.

Our approach contrasts with the usual approach (e.g., \cite{wahba90}, Section 10.2, \cite{bt04}, Section 10.2, \cite{gu13}, Section 2.4.5), which is based on a kernel
\begin{align}\label{non-pars}  h(x,x') =\sum_{\cM\in\tilde\dP(\dM)}\tau_\cM\tilde h_\cM(x_\cM,x_\cM'), \end{align}
where the $\tau_\cM$ are positive and otherwise unrestricted. 
Thus, if all main effects are present, in our approach there are $\abs{\cV}+1$ scale parameters, regardless of the number of interactions in the model, while in the usual approach, there are 
$\abs{\tilde\dP(\dM)}\ge\abs{\cV}+1$ scale parameters. Note that~(\ref{non-pars}) does not define a hierarchical interaction space unless there exist $\xi$ and $\lambda_v$ such that the factorization  $\tau_\cM=\xi\prod_{v\in\cM}\lambda_v$ holds. Kernels of of the forms discussed above have also been referred to as {\em ANOVA} kernels.


\subsection{Hierarchical interaction models} \label{sec-hier2}


We assume observations $(x_{1},y_1),\ldots,(x_{n},y_n)$, where $x_{1},\ldots,x_n\in\cX_\cV$ and $y_1,\ldots,y_n\in\mR$. Let $\dM$ be a Sperner family of $\cV$. In the sequel, we refer to such $\dM$ as the {\em maximal interactions} of a model. 
\begin{definition}[Hierarchical interaction model]
The observations $(x_{1},y_1),\ldots,(x_{n},y_n)$ satisfy the {\em hierarchical interaction model} associated with the {maximal interactions} $\dM$ if they satisfy
\begin{align}\label{regr}
	y_i = f(x_i) + \ve_i, \quad f\in\cF_\dM,  i=1,\ldots,n  ,
\end{align}
where $\ve_i\in\mR$ is a random error term. 
The $f_\cM$ with $\abs{\cM}=1$ are called {\em main effects} and the $f_\cM$ with $\abs{\cM}>1$ {\em interaction effects} of order $\abs{\cM}-1$. 
The model is called {\em saturated} if $\dM=\{\cV\}$. 
\end{definition}
\noindent For example, if $\cV=\{1,2\}$ the saturated model ($\dM=\{\{1,2\}\}$) is 
\[  y_i = f_\emptyset() + f_{1}(x_{i,1}) + f_{2}(x_{i,2}) + f_{12}(x_{i,1},x_{i,2}) + \ve_i, \quad i=1,\ldots,n  \]
where $x_i=(x_{i,1},x_{i,2})$, $f_{\{1\}}$ and $f_{\{2\}}$ are the main effects, and $f_{\{1,2\}}$ is the interaction effect. 
The non-saturated model associated with the maximal interactions $\dM=\{\{1\},\{2\}\}$ is
\[  y_i = f_\emptyset() + f_{1}(x_{i,1}) + f_{2}(x_{i,2}) + \ve_i, \quad i=1,\ldots,n  \]
Furthermore, if $\dM=\emptyset$, $y_i=\ve_i$ and if $\dM=\{\emptyset\}$, $y_i=f_\emptyset()+\ve_i$, $i=1,\ldots,n$ (the intercept only model). 

In addition to A1 and A2 above, we make the following assumption in this paper: 
\begin{enumerate}
    \item[A3] The errors in~(\ref{regr}) have a multivariate normal distribution, i.e.,
\begin{eqnarray}\label{err}
	(\varepsilon_1,\ldots,\varepsilon_n) \sim \MVN(0,\Psi^{-1}),
\end{eqnarray}
where $\Psi=(\psi_{ij})$ is an $n\times n$ positive definite precision matrix. Here, $\Psi$ is taken to be known up to a low dimensional parameter, e.g., $\Psi=\psi I_n$ ($\psi>0$, $I_n$ the $n\times n$ identity matrix), reflecting i.i.d.\ errors.

\end{enumerate}

\section{Estimating hierarchical interaction models with I-priors} \label{sec-noninf}

In this section we consider model~(\ref{regr}) subject to Assumption A3 and where $\cF_\cM$ is the RKKS with reproducing kernel~(\ref{parsim kern}). (Note that, therefore, Assumptions A1 and A2 are also satisfied.)
We define and derive the I-prior for the regression function $f$, then give the marginal likelihood and an efficient EM algorithm for estimating the scale parameters of the kernels. Note that the derivation of the I-prior {\em only} requires the previously stated model assumptions, that is, no further user input is needed.

\subsection{I-priors}

As shown in Appendix~\ref{app-fisher}, the Fisher information on $f$ is given by
\[  I[f](x,x') = \sum_{i=1}^n\sum_{j=1}^n \psi_{ij}h(x,x_i) h(x,x_j). \]
$I[f]$ is positive definite and hence induces an RKHS over $\cX$, which we denote by $\cF_n$. $\cF_n$ is a finite dimensional subspace of $\cF$, consisting of functions of the form $f(x)=\sum_{i=1}^nh(x,x_i)w_i$, with squared norm $\norm{f}_{\cF_n}^2=\sum\sum w_iw_j\psi_{ij}^{-}$, where $\psi_{ij}^-$ is the $(i,j)$th element of the error covariance matrix $\Psi^{-1}$ \citep[Lemma~2]{bergsma20}.
Note that the standard score test statistic is the $\cF_n$-norm of the score vector.



It is immediate that we can uniquely decompose any $f\in\cF$ as $f=f_n+r_n$, where $f_n\in\cF_n$ and $r_n(x_1)=\ldots=r_n(x_n)=0$. 
Such an $r_n$ is orthogonal to $\cF_n$ in $\cF$ since, by the reproducing property of $h$, $\langle f_n, r_n\rangle_\cF=\sum w_i\langle h(x_i,\cdot),r_n\rangle_\cF=\sum w_ir_n(x_i)=0$. Hence, $r_n\in\cF_n^\perp$, where $\cF_n^\perp$ is the orthogonal complement of $\cF_n$ in $\cF$. Since additionally the likelihood of $f$ clearly does not depend on $r_n$, we can conclude the data contain no information on $r_n$ (alternatively formulated, the Fisher information for $r_n$ is zero, see also \cite[Lemma~3]{bergsma20}). Therefore, any `objective' prior, i.e., a prior not determined by subjective beliefs, should have its support restricted to $\cF_n$, and we can only estimate $r_n$ by a `best guess'; in this paper, we set $r_n=0$. 

We define the I-prior as a maximum entropy prior as follows. 
Let $\nu$ be volume measure induced by $\norm{\cdot}_{\cF_n}$ (since $\norm{\cdot}_{\cF_n}$ does not depend on $f$, this measure is flat). 
The entropy of a prior $\pi$ over $\cF_n$ relative to $\nu$ is
\[  {\cal E}(\pi) = - \int_{\cF_n}\pi(f)\log\pi(f)\nu(\diff f). \]
We define the I-prior for $f$ as the prior $\pi$ maximizing entropy subject to the constraint
\[  \E_{g\sim \pi}\norm{g}_{\cF_n}^2= \mbox{constant}. \]
Variational calculus shows that an I-prior for $f$ is the Gaussian variable with mean  and covariance kernel proportional to the Fisher information on $f$, i.e., 
\[  \cov_{\pi}(f(x),f(x')) \propto \sum_{i=1}^n\sum_{j=1}^n \psi_{ij}h(x,x_i)h(x',x_j). \]
Equivalently, under the I-prior, $f$ can be written in the form
\begin{eqnarray}  \label{repr2}
	f(x) \propto \sum_{i=1}^n h(x,x_i)w_i, \hspace{10mm}(w_1,\ldots,w_n)\sim\MVN(0,\Psi).
\end{eqnarray}
In this paper, we use the kernel~(\ref{parsim kern}) which has an overall scale parameters $\tau$, so with this kernel we can replace proportionality in the above two equations with equality. 

In \citet{bergsma20}, it is shown that $g$-priors \citep{zellner86} can be viewed as a special case of I-priors, if the regression function is linear over $\cX=\mR^p$ equipped with Mahalanobis distance. 
Furthermore, he showed that if $\cX=\mR$ and the regression function is in the Brownian motion RKHS (the RKHS whose reproducing kernel is the Brownian motion covariance kernel), the posterior mean of the regression function under the I-prior is closely related to cubic spline smoothers.

\subsection{Marginal likelihood and posterior distribution of parameter estimates}\label{sec-marg}

Denote $\by=(y_1,\ldots,y_n)^\top$, $\bff=(f(x_1),\ldots,f(x_n))^\top$, $\bff_0=(f_0(x_1),\ldots,f_0(x_n))^\top$, $\bvarepsilon=(\varepsilon_1,\ldots,\varepsilon_n)^\top$, and $\bw=(w_1,\ldots,w_n)^\top$. 
Then~(\ref{regr}) implies $\by = \bff + \bvarepsilon$.
Let $\bH_{\blambda}$ be the $n\times n$ matrix with $(i,j)$th coordinate $h_{\blambda}(x_i,x_j)$, where $h_{\blambda}$ is the reproducing kernel of $\cF_\dM$  with scale parameter vector $\blambda\in\mR^{\abs{\cV}}$; that is, $h_\blambda$ is given by the right-hand side of~(\ref{parsim kern}).

Under the I-prior, $\bff \sim \MVN(\bff_0,\bH_{\blambda}\bPsi \bH_{\blambda})$, so the
marginal distribution of $\by$ is 
\begin{eqnarray}\label{ydist}
	\by \sim \MVN(\bff_0,\bV_{\by}),
\end{eqnarray}
where the marginal covariance is given as
\begin{align}\label{vydef}
	\bV_{\by} = \bH_{\blambda}\bPsi\bH_{\blambda}+\bPsi^{-1} .
\end{align}
Thus, the marginal log likelihood of $(\blambda,\Psi)$ is
\begin{equation}\label{marglik}  
	l(\blambda,\bPsi|\by) =  -\frac{n}{2}\log(2\uppi) - \frac12\log|\bV_\by| - \frac12(\by-\bff_0)^\top \bV_\by^{-1}(\by-\bff_0). 
\end{equation}
The maximum likelihood (ML) estimate $(\hat\blambda,\hat\bPsi)$ of $(\blambda,\bPsi)$ maximizes $l(\blambda,\bPsi|\by)$, and can be obtained using the EM algorithm described below. 
Having obtained $(\hat\blambda,\hat\bPsi)$, $f$ can be estimated by its posterior distribution under the I-prior:
\begin{lemma}\label{lem-post}
	The posterior distribution of $f$ in~(\ref{regr}) subject to~(\ref{err}) given $\by$ under the I-prior $\pi$ is Gaussian with mean given by
	\[  \E_\pi\big[f(x) \,|\, \by\big] = f_0(x) + \sum_{i=1}^n h(x,x_i)\hat w_i,  \]
	where
	\begin{align}\label{tildew}
		\hat\bw = \bPsi\bH_\blambda\bV_\by^{-1}(\by-\bff_0) ,
	\end{align}
	and covariance kernel given by 
	\begin{align*}
		\cov_{\pi}\big(f(x),f(x')|y_1,\ldots,y_n\big) = \sum_{i=1}^n\sum_{j=1}^n h(x,x_i)h(x',x_j)(\bV_{\by}^{-1})_{ij}.
	\end{align*}
\end{lemma}
\noindent The proof uses standard methods and is given in \citet{bergsma20}. 


\subsection{Estimation of interaction models using the EM algorithm} \label{sec-em}


We now describe the EM algorithm for estimating the scale parameter $\blambda$ of $\bH_{\blambda}$, as well as parameters of the precision matrix $\bPsi$. For estimating these parameters, EM turns out to be particularly efficient. 
The E-step is in closed form and the M-step is computationally very easy to carry out, involving the solution of several polynomials in several variables.
Kernels may have additional hyperparameters, e.g., the Hurst coefficient for the fractional Brownian motion kernel, but our EM algorithm is computationally intensive for estimating these. 


With $g$ denoting the density function related to its argument and using~(\ref{repr2}) with $\tau=1$, the complete data log likelihood is
\begin{align*}
	l(\blambda,\bPsi|\by,\bw)
	&= \log g(\by|\bw,\blambda,\bPsi) + \log g(\bw|\bPsi) \\
	&= c + \cancel{\frac12\log|\Psi|} - \frac12(\by-\bff)^\top\bPsi(\by-\bff) - \cancel{\frac12\log|\Psi|} - \frac12\bw^\top\bPsi^{-1}\bw   \\
	&= c - \frac12(\by-\bff)^\top\bPsi(\by-\bff) - \frac12\bw^\top\bPsi^{-1}\bw   \\
	&= c - \frac12(\by-\bff_0-\bH_{\blambda}\bw)^\top\bPsi(\by-\bff_0-\bH_{\blambda}\bw) - \frac12\bw^\top\bPsi^{-1}\bw   \\
	&= c - \frac12(\by-\bff_0)^\top\bPsi(\by-\bff_0) - \frac12\bw^\top\bH_{\blambda}^\top\bPsi\bH_{\blambda}\bw + (\by-\bff_0)^\top\bPsi\bH_{\blambda}\bw - \frac12\bw^\top\bPsi^{-1}\bw \\
	&= c - \frac12(\by-\bff_0)^\top\bPsi(\by-\bff_0) - \frac12\bw^\top\bV_\by\bw + (\by-\bff_0)^\top\bPsi\bH_{\blambda}\bw\\
	&= c - \frac12(\by-\bff_0)^\top\bPsi(\by-\bff_0) - \frac12\tr\big[\bV_\by\bw\bw^\top\big] + (\by-\bff_0)^\top\bPsi\bH_{\blambda}\bw,
\end{align*}
where $c$ is a constant. Write
\[  \tilde\bW = E\big(\bw\bw^\top \big |\, \by,\blambda,\bPsi\big) = \tilde\bV_\bw+\tilde\bw\tilde\bw^\top, \]
where $\tilde\bw$ and $\tilde\bV_\bw$ are given by~(\ref{tildew}) and~(\ref{vydef}).
Let $\tilde\bw^{(0)}$ and $\tilde\bW^{(0)}$ be $\tilde\bw$ and $\tilde\bW$ with $\bPsi$ and $\blambda$ replaced by $\bPsi^{(0)}$ and $\blambda^{(0)}$.
The E-step consists of computing
\begin{eqnarray}
	Q(\blambda,\bPsi)
	&=& E\left\{ l(\blambda,\bPsi|\by,\bw) \, \Big|  \, \by,\blambda^{(0)},\bPsi^{(0)}\right\} \nonumber\\
	&=& c - \frac12(\by-\bff_0)^\top\bPsi(\by-\bff_0) - \frac12\tr\big[\bV_\by\tilde\bW^{(0)}\big] + (\by-\bff_0)^\top\bPsi\bH_{\blambda}\tilde\bw^{(0)}.\label{QI}
\end{eqnarray}

The M-step entails maximizing $Q(\blambda,\bPsi)$. We assume the global maximum can be found by differentiating, equating to zero, and solving.
Supposing $\bPsi$ but not $\bH_{\blambda}$ depends on a parameter $\psi$ and $\bH_{\blambda}$ but not $\bPsi$ depends on a parameter $\lambda$, the derivatives are given by
\begin{eqnarray*}
	\frac{\partial Q(\blambda,\bPsi)}{\partial\lambda}
	&=& - \tr\Big[\frac{\partial\bH_{\blambda}}{\partial\lambda}\bPsi\bH_{\blambda}\tilde\bW^{(0)}\Big] + (\by-\bff_0)^\top\bPsi\frac{\partial\bH_{\blambda}}{\partial\lambda}\tilde\bw^{(0)}  \\
	\frac{\partial Q(\blambda,\bPsi)}{\partial\psi}
	&=& - \frac12(\by-\bff_0)^\top\frac{\partial\bPsi}{\partial\psi}(\by-\bff_0)
	- \frac12\tr\big[\frac{\partial\bV_\by}{\partial\psi}\tilde\bW^{(0)}\big]
	+ (\by-\bff_0)^\top\frac{\partial\bPsi}{\partial\psi}\bH_{\blambda}\tilde\bw^{(0)}.
\end{eqnarray*}
For the example in Section~\ref{sec-cow}, the errors are assumed i.i.d., so $\bPsi=\psi\bI_n$ for a scalar $\psi$, so $\bH_{\blambda}=\sum_{v\in\cV}g_v(\blambda)\bH_{(v)}$, where $g_s$ is a polynomial function and $\bH_{(v)}$ is the Gram matrix for the kernel $h_v$.
Then
\begin{eqnarray*}
	\frac{\partial\bPsi}{\partial\psi} = \bI_n   \hspace{10mm} \frac{\partial\bH_{\blambda}}{\partial\lambda} = \sum_{v\in\cV}\frac{\partial g_v(\blambda)}{\partial\lambda}\bH_{(v)}.
\end{eqnarray*}

The partial derivatives of $Q$ set to zero can then normally be solved very quickly numerically for the purposes of this paper.
For example, if $g_s(\blambda)=\lambda_s$ (which is the case, e.g., if there are no interactions in the model) then the equations have a closed form solution. 
In general, with i.i.d.\ errors, and $\cV$ the set of (possibly multidimensional) variables involved, $2|\cV|+1$ polynomial equations in $2|\cV|+1$ unknowns need to be solved, which we found can be done very rapidly using built in solvers in R and Mathematica. The computational bottleneck is not the M step, but the E step which is $O(n^3)$. 


The EM algorithm described above is computationally unattractive for estimating kernel hyperparameters. We note however that reasonably flexible analysis can be done using kernels without hyperparameters, in particular, using linear kernels and the multivariate Brownian motion kernel. The latter consists of functions possessing a directional derivative \citep{bergsma20} which may be attractive for many applications. For the data example in Section~\ref{sec-cow}, we used the fractional Brownian motion RKHS and estimated the Hurst coefficient using a bisection algorithm for the profile marginal likelihood. 

A difficulty with estimation is that the marginal likelihood for the scale parameters may have multiple local maxima, increasing with the number of scale parameters. Our parsimonious approach which has no additional scale parameters for interaction effects reduces this difficulty, but still typically multiple starting points need to be tried with the EM algorithm. Further research is needed to determine good starting values. 

\section{Simulation study on model selection}\label{sec-sim}



The goal of this simulation study is to compare how well different methods are at identifying the presence or absence of main effects, two-way interactions and the single three-way interaction in a (hierarchical) linear model. Models with linear regression functions are a useful way to evaluate the general I-prior methodology because for this case there are well-developed methods to compare with, such as the lasso, spike-and-slab priors, and $g$-priors. 

Data pairs $\{(x_i, y_i)\}_{i=1}^n$, where $x_i=(x_{i1},x_{i2},x_{i3})\in\mathbb{R}^3$, were simulated for $i=1,\dots,n=100$ according to the following generating process:
\vspace*{-\baselineskip}
\begin{equation}
    y_i = \beta_1 x_{i1} + \beta_2 x_{i2} + \beta_3 x_{i3} + \beta_{12}
    x_{i1}x_{i2} + \beta_{13} x_{i1}x_{i3} + \beta_{23} x_{i2}x_{i3} +
    \beta_{123} x_{i1}x_{i2}x_{i3} + \epsilon_i.
\end{equation}
Note here that the errors $\epsilon_i$ were independently and identically generated from the $\N(0,\sigma^2$) distribution, with $\sigma=3$, while the independent variables $x_i$ were generated according to
\begin{equation}
    \begin{pmatrix}
    x_{i1}\\
    x_{i2}\\
    x_{i3}
    \end{pmatrix}
    \sim \MVN
    \left(
    \begin{pmatrix}
    0\\
    0\\
    0
    \end{pmatrix}
    ,
    \begin{pmatrix}
    1 &\rho &\rho \\
    \rho & 1 & \rho\\
    \rho & \rho & 1
    \end{pmatrix}
    \right).
\end{equation} \vspace{-0.5em}
It is of interest to see the performance of each method when dealing with both uncorrelated data ($\rho=0$) and correlated data ($\rho=0.5$).

The coefficients $\beta_1,\beta_2,\dots,\beta_{123}$ were varied to produce data from eight different models, as per Table \ref{tab:betas}  below.
Each of the eight models represent an interesting realization of the hierarchical linear model. 
For example, model 4 (1110000) is the model corresponding to the presence of main effects only, while model 5 (1110100) is the model corresponding to the presence of main effects with a single two-way interaction between $x_1$ and $x_3$. 
While we could have also looked at the two other possibilities for the single two-way interaction case, namely 1111000 and 1110010, there really wouldn't be any reason to prefer one over any other, and we would therefore expect similar results in each case.
This, in principle, ruled out the other possible models to look at, reducing the number of interesting models to only the eight listed in Table \ref{tab:betas}.

\begin{table}[ht] 
\centering
\caption{Coefficients for the different models of the data generating process. The model code here is a 7-digit binary string identifying the inclusion (1) or exclusion (0) of the corresponding effect in that digit's position. }
\label{tab:betas}
\begin{tabular}{@{}lcccccccr@{}}
\toprule
No.  & $\beta_1$ & $\beta_2$ & $\beta_3$ & $\beta_{12}$ & $\beta_{13}$ & $\beta_{23}$ & $\beta_{123}$ & Model code    \\ \midrule
1 & 1  & 0  & 0  & 0  & 0  & 0  & 0   & 1000000 \\
2 & 1  & 1  & 0  & 0  & 0  & 0  & 0   & 1100000 \\
3 & 1  & 1  & 0  & 0.5  & 0  & 0  & 0   & 1101000 \\
4 & 1  & 1  & 1  & 0  & 0  & 0  & 0   & 1110000 \\
5 & 1  & 1  & 1  & 0  & 0.5  & 0  & 0   & 1110100 \\
6 & 1  & 1  & 1  & 0.5  & 0.5  & 0  & 0   & 1111100 \\
7 & 1  & 1  & 1  & 0.5  & 0.5  & 0.5  & 0   & 1111110 \\
8 & 1  & 1  & 1  & 0.5  & 0.5  & 0.5  & 0.25   & 1111111 \\ \bottomrule
\end{tabular}
\end{table}

Four competing methods were used to fit the simulated data set coming from each model, namely: 1) our I-prior method; 2) the lasso \citep{glmnet}; 3) the fully Bayesian spike-and-slab variable selection method; and 4) the $g$-prior method.
All methods were implemented in R code, which is provided as a supplementary material to this paper.
In the case of the I-prior method, all hierarchical interaction models were fitted and the model giving the highest log-likelihood was selected. For a fair comparison with the competing methods, observed covariates were standardized, hence the same scale parameter could be used for each. This left three hyperparameters to be estimated: the common scale parameter $\lambda$, the overall scale parameter $\tau$, and the error variance.
Estimation was done with maximum marginal likelihood using the R package \texttt{iprior} \citep{jamil19package}.
In the case of the lasso, the tuning parameter $\lambda$ was chosen by 10-fold cross-validation, and  the $\lambda$ value giving the smallest cross-validation error used.
The spike-and-slab approach was done using MCMC sampling, making use of the JAGS functionality \citep{plummer2003jags} from the \texttt{runjags} package \citep{denwood2016runjags}.
Uninformative priors were used: $\N(0,100)$ for each regression coefficient $\beta$, $\Gamma(0.01,0.01)$ for the error variance $\sigma^2$, and uniformly independent Bernoulli distributions for the model indicators. 
The final model selected was based on maximal model posterior probability.
For more details on this method, see for example \citet{George1993}, \citet{Ntzoufras2008}, or \citet{jamil_bayesian_2021}.
Lastly, the $g$-prior method was fitted using the \texttt{BAS} package \citep{clydeBAS}, where the estimation of the hyperparameter $g$ was done using the MLE of $g$ from the marginal likelihood within each model considered.

For each of the eight models, $B=10,000$ data replications were generated and each of the four methods' results saved for comparison.
Of interest here is the proportion of times the correct data generating process was reported by each method for each of the models listed above.
A meaningful summary of the performance of the four methods would be the geometric mean of correct proportions, as depicted in Figure \ref{fig:results} below.
Evidently, the I-prior is the superior method of the four, performing the best in all scenarios tested.

\begin{figure}[htbp]
	\centering
	\includegraphics[width=0.75\textwidth]{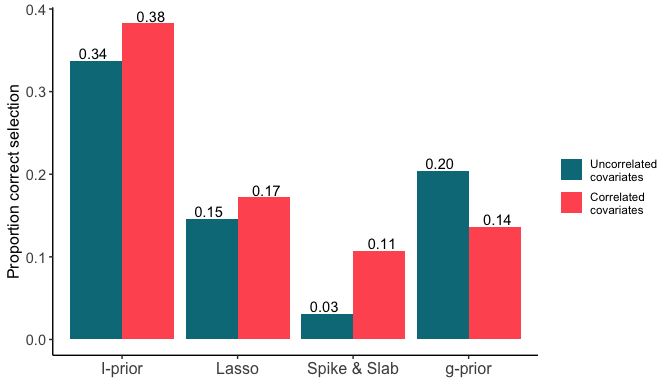}
	\caption{Geometric mean of probabilities of correct model selection across the two simulation scenarios (i.e., correlated and uncorrelated covariates). The list of models simulated from is given in Table~\ref{tab:betas}. }
 \label{fig:results}
\end{figure}

Most methods performed similarly in the presence of correlated or uncorrelated covariates. 
Perhaps this was slightly surprising, given the known fact that the lasso performs poorly under multicollinearity. 
A closer inspection of the simulation results (details of which are given in the supplementary material) reveals that the $g$-prior frequently selects the saturated model (1111111) under any circumstance involving interactions.
This error contributed to the high proportion of selection when the true model was indeed the saturated model, so thus inflating the results slightly in its favour.
The worst performing method was the independent spike-and-slab one in a fully Bayesian setting.
This method had consistent difficulty in selecting the correct interaction effect, but generally no issues when only main effects were present.


\section{Application: a functional response model} \label{sec-cow}



\begin{figure}[htbp]
	\centering
        \includegraphics[width=0.8\textwidth]{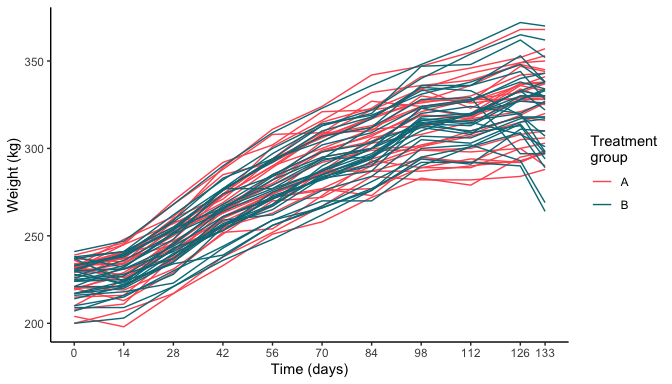}
	\caption{
 Growth curves of cows. The two colours represent the two treatments received. The question of interest is whether and how treatment affects growth.}\label{fig-cow}
\end{figure}

We consider a balanced longitudinal data set consisting of repeated measurements of weights of 60 cows, 30 of which are randomly assigned to a treatment group $A$ and 30 to a treatment group $B$. Weight was measured 11 times over a 133-day period, at two-week intervals, except for the last measurement, which was taken one week after the preceding measurement. In Figure~\ref{fig-cow} a sample of growth curves is shown.

A common approach to analyze a longitudinal data set such as this one is to assume that the observed growth curves are realizations of a Gaussian process with unknown covariance function, and to estimate the latter.
For example, \citet{kenward87} assumed a so-called ante-dependence structure of order $k$, which assumes an observation depends on the previous $k$ observations, but given these is independent of any preceding observations.
Various other process families have been considered (\cite{nz00,pourahmadi00,pm03,zlt14}; see the latter for an overview).

In contrast, the I-prior methodology relies on specification of a class of functions for the growth curves. Here, we assume that up to an intercept, the growth curves lie in a fractional Brown motion (FBM) RKHS. The FBM RKHS on $\mR$ with Hurst coefficient $\gamma$ is defined by the reproducing kernel
\[   h_\gamma(t,t') = \frac12\big(\abs{t}^{2\gamma} + \abs{t'}^{2\gamma} - \abs{t-t'}^{2\gamma}\big). \]
The RKHS consists of functions possessing a square integrable fractional derivative of order $2\gamma$ \citep{picard11}.

The weight of cow $i$ ($i=1,\ldots,m$) at time $t$ ($t=t_1,\ldots,t_T$) is denoted $y_{it}$. There are three covariates to explain weight: (i) a nominal categorical covariate {\em treatment} (treatment group of cow $i$ denoted $x_i\in\{A,B\}$); a nominal categorical covariate {\em cow index}; and a continuous covariate {\em time}. 
The regression model we consider to explain weight is of the form
\[  y_{it} = f(i,x_i,t) + \varepsilon_{it}, \]
where $f$ lies in a hierarchical interaction space with r.k.\ of the form~(\ref{parsim kern}) and the $\ve_{it}$ are assumed to be independent $\N(0,\sigma^2)$ variables. Examples of possible decompositions of $f$ are of the form~(\ref{sat3}) or~(\ref{nonsat}) (there are 20 possible hierarchical decompositions, see Section~\ref{sec-hier1}).

We next need to specify the main effect kernels for each covariate. We assume the main effects of the nominal categorical covariates treatment and cow index lie in the identity RKHS, defined as the unique RKHS with reproducing kernel the identity kernel, defined by $h(i,i')=\delta_{i,i'}$ where $\delta$ is Kronecker delta. Effectively, this is the same as coding nominal outcomes using dummy vectors and assigning Euclidean distance. The main effect of time is assumed to be in the FBM RKHS. We estimate the Hurst coefficient using the maximum marginal likelihood method.




For ease of interpretation, we reformulate the scalar-response regression model as a multivariate (or functional) response model, where the response is the growth curve of a cow. Denote the growth curve of cow $i$ by $\by_i=(y_{t_1},\ldots,y_{t_T})$, let $\bff(i,x_i)=(f(i,x_i,t_1),\ldots,f(i,x_i,t_T))$ and let $\bm{\ve}_i=(\ve_{it_1},\ldots,\ve_{it_T})$. 
Then
\[  \by_{i} = \bff(i,x_i) + \bm{\ve}_{i}, \]
Under the saturated model,
\[  \bff(i,x_i) = \balpha + \bff_1(i) + \bff_2(x_i) + \bff_{12}(i,x_i) . \]
The hierarchical models (with no restrictions on the time effect) are indexed by their maximal effects (i.e., Sperner families) $\emptyset$, $\{\emptyset\}$, $\{\{C\}\}$, $\{\{X\}\}$, $\{\{C\},\{X\}\}$, and $\{\{C,X\}\}$. We denote the latter four by the shorthands $\{C\}$, $\{X\}$, $\{C,X\}$ and $\{CX\}$. 
Under model $\emptyset$, $\bff(i,x_i)=0$ and we do not consider this further. Under model $\{\emptyset\}$, $\bff(i,x_i)=\alpha$, i.e., growth curves are the same for each cow and treatment; model $\{C\}$ asserts different growth curves for different cows, but no treatment effect; 
model $\{X\}$ asserts an effect of treatment on growth but otherwise identical growth for all cows; model $\{C,X\}$ asserts both treatment and cow index affect growth; the saturated model $\{CX\}$ additionally asserts that the effect of treatment is different for different cows. From inspection of the curves in Figure~\ref{fig-cow}, it seems clear any well-fitting model should include an effect of $C$. We are not aware of the fit of the interaction model $\{CX\}$ having been studied in the literature. 

Table~\ref{tbl-cowfit} gives a summary of the model fits. A test of model $\{C,X\}$ against $\{C\}$ is a conditional independence test (between $X$ and growth given $C$) and in general not straightforward \citep{sp20,berrett19,bergsma04a}. In the present case, $X$ (treatment) is independent of $C$ (cow index), and hence the distribution of the maximum marginal likelihood under model $\{C,X\}$ assuming model $\{C\}$ is true can be simulated using the permutation method, that is, by randomly permuting the treatment labels of cows a large number times. This yielded a simulated $p$-value of $0.000$ based on 1000 random permutations, i.e., there is strong evidence for a treatment effect. 

Our parsimonious kernel specification~(\ref{parsim kern}) makes it easy to compare models with the same main effects and different interaction effect, namely by selecting the model with the highest log-likelihood. Hence we select model $\{C,X\}$ over $\{CX\}$, i.e., the data suggest that treatment affects different cows identically.

\begin{table}[t]
 	\caption{Goodness of fit for cow data. The model consists of the highest order effects on the growth curve, e.g., model $\{C,X\}$ means the growth curve depends on cow ($C$) and treatment ($X$), and there is no interaction meaning that the treatment effect is the same for all cows.}
 	\label{tbl-cowfit}
	\centering
		\begin{tabular}{llllll}\toprule
			Model      & Log-likelihood & Error standard  & Number of $\lambda$ & AIC & BIC \\
			&                & deviation       & parameters   &  \\ \midrule
			   $\emptyset$     & $-2788.8$  & 16.3 & 1  & 5583.5 & 5597.0 \\
			   $\{X\}$    & $-2788.7$  & 16.3 & 2  & 5585.5.6 & 5603.5 \\
			   $\{C\}$    & $-2253.2$  & 0.2  & 2  & 4514.4 & 4532.4 \\
			   $\{C,X\}$  & $-2231.1$  & 0.1  & 3  & 4472.3 & 4494.7 \\
			   $\{CX\}$   & $-2232.8$  & 0.2  & 3  & 4475.6 & 4498.0 \\ \bottomrule
		\end{tabular}
\end{table}

\section{Remarks}\label{sec-disc}


The I-prior methodology allows the estimation of a regression function $f$ of a model $y_i=f(x_i)+\ve_i$ ($i=1,\ldots,n$), where $f$ lies in an RKKS $\cF$ with kernel of the form~(\ref{parsim kern}) and the errors are normal (e.g., i.i.d.\ normal or AR(p) errors). 
The I-prior methodology has the following properties:
\begin{enumerate}
	\item Since the support of the I-prior is contained in $\cF$, the complete class theorem implies the posterior mean of $f$ under the I-prior is admissible under a broad range of loss functions.
	\item The I-prior methodology requires minimal assumptions on the regression, in particular, an RKKS $\cF$ and an error distribution.
	\item An EM algorithm with simple E and M steps for finding the maximum likelihood estimators of the scale parameters $\lambda_v$ of the kernel~(\ref{parsim kern}) and of the parameters of the error distribution is available (see Section~\ref{sec-em}). 	
\end{enumerate}
The main alternatives to I-prior modelling are Gaussian process (GP) regression and Tikhonov regularization. 
The first property gives I-prior estimators an advantage over Tikhonov regularizers, which are inadmissible with respect to squared error loss \citep{cp19}. 
The second property gives the I-prior methodology an advantage over GP regression, which in addition to a metric over $\cF$ requires the user to specify a prior. 
The third property gives the I-prior methodology an advantage over both Tikhonov regularization and GP regression, for which no simple and generally applicable algorithms are available to estimate scale parameters. 
We should note that to estimate {\em kernel hyperparameters}, such as the variance parameter of a squared exponential kernel, the I-prior methodology does not have an advantage over competing methods.

\bigskip
\begin{center}
{\large\bf SUPPLEMENTARY MATERIAL}
\end{center}

\begin{description}

\item[Code for simulation:] The R codes and detailed simulation results are provided here: \url{https://github.com/haziqj/iprior-interaction}. (website) 

\end{description}



\appendix

\section{Fisher information on the regression function}\label{app-fisher}

Consider model~(\ref{regr}) subject to assumptions A1 and A2. By Proposition~\ref{prop-kern}.1, $f(x)=\langle f,h(x,\cdot)\rangle_{\cF_\dM}$, where $h$ is the r.k.\ of $\cF_\dM$ given by Proposition~\ref{prop-kern}.
The log-likelihood of $f$ in~(\ref{regr}) then is
\[   l(f) = C - \frac12\sum_{i=1}^n\sum_{j=1}^n \psi_{ij}\big(y_i-\langle f,h(x_i,\cdot)\rangle_{\cF_\dM}\big)\big(y_j-\langle f,h(x_j,\cdot)\rangle_{\cF_\dM}\big), \]
for a constant $C$. Since $\nabla_f\langle f,h(x_i,\cdot)\rangle_{\cF_\dM}=h(x_i,\cdot)$, the score function is the gradient of $l$, given as
\[   \nabla l(f) = \sum_{i=1}^n\sum_{j=1}^n \psi_{ij}\big(y_i - \langle f,h(x_i,\cdot)\rangle_{\cF_\dM}\big)h(x_j,\cdot) .  \]
The Fisher information on $f$ then is
\[   I[f] = - \nabla^2 l(f) = \sum_{i=1}^n\sum_{j=1}^n \psi_{ij}h(x_i,\cdot)\otimes h(x_j,\cdot) .  \]


\bibliography{stats}

\end{document}